\documentclass[11pt,reqno,a4paper]{amsart}

\usepackage{etoolbox}
\usepackage{amsmath}
\usepackage{enumerate}
\usepackage{amssymb}
\usepackage{amscd}
\usepackage{amsthm}
\usepackage{amsfonts}
\usepackage{graphicx}

\patchcmd{\subsection}{-.5em}{.5em}{}{}
\patchcmd{\subsubsection}{-.5em}{.5em}{}{}

\usepackage{enumitem}

\usepackage[scaled]{helvet} 
\usepackage{courier} 
\usepackage{eulervm}
\normalfont

\usepackage{hyperref}

\numberwithin{equation}{section}

\newcommand{\Aut}{\operatorname{Aut}}

\newcommand{\bN}{\mathbb{N}}

\newcommand{\bR}{\mathbb{R}}

\newcommand{\bT}{\mathbb{T}}

\newcommand{\bZ}{\mathbb{Z}}

\newcommand{\ra}{\rightarrow}

\newcommand{\qor}{\quad \textrm{or} \quad}
\newcommand{\qand}{\quad \textrm{and} \quad}

\newcommand\subsetsim{\mathrel{%
\ooalign{\raise0.2ex\hbox{$\subset$}\cr\hidewidth\raise-0.8ex\hbox{\scalebox{0.9}{$\sim$}}\hidewidth\cr}}}
\newcommand{\eps}{\varepsilon}

\DeclareMathOperator{\Stab}{Stab}

\theoremstyle{theorem}
\newtheorem{theorem}{Theorem}[section]

\newtheorem{proposition}[theorem]{Proposition}
\newtheorem{lemma}[theorem]{Lemma}

\theoremstyle{definition}
\newtheorem{definition}[theorem]{Definition}
\newtheorem{remark}[theorem]{Remark}

\begin{document}

\title{Small product sets in compact groups}

\author{Michael Bj\"orklund}

\address{Department of Mathematics, Chalmers, Gothenburg, Sweden}
\email{micbjo@chalmers.se}



\date{}


\maketitle

\begin{abstract}
We show in this paper that a sub-critical pair $(A,B)$ of sufficiently "spread-out" Borel sets in a compact and second countable group $K$ with an \emph{abelian} identity component, must reduce to a Sturmian pair 
in either $\bT$ or $\bT \rtimes \{-1,1\}$. This extends a classical result of Kneser.
\end{abstract}

\section{The main result}

Let $K$ be a compact and second countable Hausdorff group with Haar probability measure $m_K$. 
Given two subsets $A, B \subset K$, we define their \textbf{product set} $AB$ by
\[
AB = \big\{ ab \, : \, a \in A, \enskip b \in B \big\}.
\]
If $A$ and $B$ are both Borel measurable subsets, then $AB$ is always measurable with respect to the \emph{$m_K$-completion} of the Borel $\sigma$-algebra on $K$, but it might fail to be Borel measurable.
This technical point will not play a major role in this paper.  \\

Before we proceed, we stress that all groups that we shall consider in this paper are assumed to be Hausdorff.

\begin{definition}[Critical and sub-critical pairs]
Suppose that $A, B \subset K$ are Borel sets. We say that $(A,B)$ is \textbf{critical} if
\[
m_K(A), m_K(B) > 0 \qand m_K(AB) <  m_K(A) + m_K(B), 
\]
and \textbf{sub-critical} if 
\[
m_K(A), m_K(B) > 0 \qand m_K(AB) =  m_K(A) + m_K(B) < 1.
\]
We warn the reader that the opposite nomenclature concerning these types of sets are sometimes 
adopted in the literature. 
\end{definition}

\begin{definition}[Reduction]
Let $M$ be a factor group of $K$ and suppose that $I, J \subset M$ are Borel sets. We say 
that $(A,B)$ \textbf{reduces} to $(I,J)$, if
\[
A \subset q^{-1}(I) \qand B \subset q^{-1}(J) \qand m_K(AB) = m_M(IJ),
\]
where $q : K \ra M$ denotes the canonical quotient map.
\end{definition}

We 
denote by $\bT$ the one-dimensional torus group $\bR/\bZ$, and by $\bT \rtimes \{-1,1\}$ the (non-abelian) 
semi-direct product of $\bT$ with the multiplicative group $\{-1,1\}$.

\begin{definition}[Sturmian pair]
Let $I, J \subset \bT$ be closed and symmetric intervals and assume that $m_{\bT}(I) + m_{\bT}(J) < 1$.
Let $M$ denote either $\bT$
or $\bT \rtimes \{-1,1\}$. We say that a pair $(A,B)$ of Borel sets in $M$ is \textbf{Sturmian}
if there exist $a, b \in M$ such that either
\[
(A,B) = (aI,Jb) \qor (A,B) = \big(a(I \rtimes \{-1,1\}),(J \rtimes \{-1,1\})b\big).
\]
One readily verifies that every Sturmian pair is sub-critical in $M$.

\end{definition}
In their very influential papers \cite{K64}, \cite{K56}, Kemperman and Kneser established:

\begin{theorem}[Theorem 1, \cite{K64}, and Satz 4, \cite{K56}]
\label{knes}
Suppose that $A, B \subset K$ are Borel sets. If 
\begin{itemize}
\item $(A,B)$ is critical, then there exists an open normal subgroup $U \triangleleft K$ such that $ABU = AB$. 
In particular, $AB$ is clopen, and $(A,B)$ reduces to a pair in a \emph{finite} factor group of $K$. 
\item $K$ is \emph{connected} and \emph{abelian}, and $(A,B)$ is sub-critical, then it
reduces to a Sturmian pair in $\bT$.
\end{itemize}
\end{theorem}

\begin{remark}
\label{AIBJ}
The first assertion effectively reduces the study of critical pairs in $K$ to the study of critical
pairs in finite groups, where powerful results of Kemperman \cite{K60}, Vosper \cite{Vo} and DeVos
\cite{DV} are available. 
\end{remark}

Recently, Griesmer \cite{G11} was able to further advance the description of sub-critical pairs in compact 
\emph{abelian} groups. Motivated by his work, and by some recent applications in ergodic theory 
related to actions of countable and discrete amenable groups developed by the author and A. Fish in \cite{BF6}, we turn in this paper 
our attention to sub-critical pairs in compact and second countable groups with an \emph{abelian} 
identity component. The relevance of this class of groups in the setting of \cite{BF6} stems 
from the observation that any compact group which contains a dense countable \emph{amenable} subgroup must have an abelian 
identity component. For a proof of this observation, we refer the reader to 
the appendix. 

\begin{definition}[Spread-out set]
A Borel subset $B \subset K$ is \textbf{spread-out} if every conull subset of $B$
projects onto every finite factor group of $K$. If $A, B \subset K$ are Borel sets, 
we say that the pair $(A,B)$ is spread-out if both $A$ and $B$ are spread-out.
\end{definition}

\begin{remark}
We stress that $K$ is always a factor of itself, and thus a \emph{proper} subset of a finite group can never 
be spread-out. Furthermore, if $(A,B)$ is a critical pair 
in $K$ and $m_K(AB) < 1$, then neither $A$ nor $B$ is spread-out by Theorem \ref{knes}, since they project onto 
the \emph{proper} subsets $AU$ and $BU$ respectively in the finite factor group $G/U$.
\end{remark}

Our main result can now be formulated as follows:

\begin{theorem}
\label{main}
Let $K$ be a compact and second countable group with an \emph{abelian} identity component. Let 
$A, B \subset K$ be Borel sets and suppose that $(A,B)$ is spread-out and sub-critical. Then 
$(A,B)$ reduces to a Sturmian pair.
\end{theorem}

\begin{remark}
We stress that this is a much weaker structural result than the previously mentioned results of Kneser, Kemperman, DeVos and Griesmer. For instance, our assumptions can never be satisfied when the identity component of $K$ 
is trivial, e.g. if $K$ is a finite group. Indeed, first note that if $A$ is spread-out in $K$ and $U \triangleleft K$ is an open normal
subgroup, then we must have $A \cap xU \neq \emptyset$ for every $x \in K/U$, since otherwise $A$ would not project onto the finite factor group $K/U$. If the identity component of $K$ is trivial, then open normal subgroups form a neighborhood basis for the identity in $K$, and thus every spread-out set in $K$ must be dense. It is not hard to see that the product set of any dense subset of $K$ and a Borel subset of $K$ of positive measure is conull. We conclude that if $A \subset K$ is spread-out and $B \subset K$ is any Borel set of positive measure, then $AB$ is conull, so in particular, $(A,B)$ cannot be sub-critical. 
\end{remark}

\section{An outline of Theorem \ref{main}}

The aim of this section is to reduce the proof of Theorem \ref{main} to two main propositions which 
will be proven in Section \ref{proofrel} and Section \ref{proofrelmain} below. \\

Let $K$ be a compact and second countable group with Haar probability measure $m_K$ and 
identity component $N$. We recall that $N$ is a closed normal subgroup of $K$, but we stress 
that it does not need to split $K$ into a semi-direct product. However, this is not far from the truth, as the 
following result by Lee shows:

\begin{theorem}[\cite{DHL}]
\label{dhl}
There exists a closed and totally disconnected subgroup $L < K$ such that $NL = K$. In particular, the semi-direct product 
group $G = N \rtimes L$, where $L$ acts on $N$ by conjugation, factors onto $K$, and we denote
by $p$ the canonical quotient map from $G$ onto $K$. 
\end{theorem}

\subsection{Sub-criticality with respect to a subgroup}

Let $G, N, L$ and $p$ be as in Theorem \ref{dhl}. We shall view $N$ and $L$ as closed subgroups of $G$. 
If $U < G$ is a closed subgroup, we consider the Haar probability measure $m_U$ on $U$ as a Borel probability measure on $G$ which is supported on $U$. 

\begin{definition}[Sub-critical with respect to a subgroup]
\label{defsubcrit}
Let $G$ be a compact and second countable group, and let $U \triangleleft G$ be a closed normal subgroup. 
We say that a pair $(A,B)$ of Borel sets in $G$ is \textbf{sub-critical with respect to $U$} if it is a sub-critical 
pair in $G$ and there are conull Borel sets $X \subset G$ and $Y \subset X \times X$ such that
\begin{equation}
\label{1}
m_U(s^{-1}A \cap U) = m_G(A) \qand m_U(Bt^{-1} \cap U) = m_G(B), 
\end{equation}
for all $s,t \in X$ and
\begin{equation}
\label{2}
m_U\big(\big(s^{-1}A \cap U\big)\big(Bt^{-1} \cap U\big)\big)
= 
m_U\big(s^{-1}A \cap U\big) + m_U\big(Bt^{-1} \cap U\big) 
\end{equation}
and 
\begin{equation}
\label{3}
m_U\big(\big(s^{-1}A \cap U\big)\big(Bt^{-1} \cap U\big)\big) 
= 
m_U\big(s^{-1}ABt^{-1} \cap U\big)
\end{equation}
for all $(s,t) \in Y$. 
\end{definition}

\begin{remark}
Since $m_G(A) + m_G(B) < 1$, we have
\[
m_U\big(\big(s^{-1}A \cap U\big)\big(Bt^{-1} \cap U\big)\big)
= 
m_U\big(s^{-1}A \cap U\big) + m_U\big(Bt^{-1} \cap U\big) 
=
m_G(A) + m_G(B) < 1.
\]
In particular, $(s^{-1}A \cap U, Bt^{-1} \cap U)$ is a sub-critical pair in $U$ for every 
$(s,t) \in Y$, which at least partially motivates the terminology "sub-critical with respect to $U$".
\end{remark}

The proof of Theorem \ref{main} breaks into two parts, guided by the following two propositions which will be established in Section \ref{proofrel} and Section \ref{proofrelmain} respectively.

\begin{proposition}
\label{relativization}
Let $A, B \subset G$ be Borel sets and suppose that $(A,B)$ is spread-out and sub-critical. 
Then $(A,B)$ is sub-critical with respect to $N$.
\end{proposition}

\begin{proposition}
\label{relativemain}
Let $A, B \subset G$ be Borel sets and suppose that $N$ is non-trivial and \emph{abelian}. If $(A,B)$ is spread-out and sub-critical with respect to $N$, then there are conull Borel sets $A' \subset A$ and $B' \subset B$ such that $(A',B')$ reduces to a Sturmian pair. 
\end{proposition}

\subsection{Proof of Theorem \ref{main}}

Let $K$ be a compact and second countable group with a non-trivial \emph{abelian} identity component $N$. Let $A_o, B_o \subset K$ be Borel sets of positive Haar measures, and suppose that $(A_o,B_o)$ is spread-out and sub-critical. 
Let $G, N, L$ and $p$ be as in Theorem \ref{dhl}. We define the Borel sets $A, B \subset G$ by
\begin{equation}
\label{defAB}
A = p^{-1}(A_o) \qand B = p^{-1}(B_o).
\end{equation}
One readily verifies that $(A,B)$ is sub-critical in $G$. 
\begin{lemma}
\label{Aisspreadout}
$(A,B)$ is spread-out in $G$.
\end{lemma}
\begin{proof}
We argue by contradiction (the proof that $B$ is spread-out is completely analogous). Suppose that we can find a conull subset $A' \subset p^{-1}(A_o)$ and a finite factor 
group $Q$ of $G$ such that $r(A') \neq Q$, where $r : G \ra Q$ denotes the canonical factor map. Since
\[
m_G(A' t \cap p^{-1}(A_o)) = m_G(p^{-1}(A_o)), \quad \textrm{for all $t \in \ker p$},
\]
and $A' \subset A' \ker r$, we have
\[
m_G((A' \ker r)t \cap p^{-1}(A_o)) = m_G(p^{-1}(A_o)), \quad \textrm{for all $t \in (\ker p)\ker r$}.
\]
Since $Q = G/\ker r$ is a finite group, the (finite) intersection
\[
A'' = \bigcap_{t \in (\ker p)\ker r} (A' \ker r)t \cap p^{-1}(A_o)
\]
is a measurable conull subset of $p^{-1}(A_o)$ which is invariant under $\ker p$ on the right hand side. Since $A'' \subset A' \ker r$, we have $r(A'') \subset r(A')$, and since $r(A') \neq Q$ and $r(A'') r(\ker p) = r(A'')$, we have
\[
r(\ker p) \neq Q \qand r(A'')/r(\ker p) \neq Q/r(\ker p).
\] 
In particular, if we let $s : K \ra Q/r(\ker p)$ denote the canonical quotient map, then $s(p(A'')) = r(A'')/r(\ker p)$ is a proper subset of the finite factor group $Q/r(\ker p)$ of $K$. Since $A'' \subset p^{-1}(A_o)$ is conull, we see that $p(A'') \subset A_o$ is a conull subset. Since $p(A'')$ does not project onto $Q/r(\ker p)$ under $s$, we see that 
the set $A_o$ is not spread-out in $K$, which contradicts our assumption about $A_o$.
\end{proof}

The assumptions of Proposition \ref{relativization} are now satisfied and we conclude that $(A, B)$ is sub-critical 
with respect to $N$. By Proposition \ref{relativemain} there are conull subsets $A' \subset A$ and $B' \subset B$ such that $(A',B')$ reduces to a Sturmian pair in either $M = \bT$ or $M = \bT \rtimes \{-1,1\}$. We recall that this means that there exists a surjective continuous homomorphism $q : G \ra M$ 
such that
\begin{equation}
\label{defCD}
A' \subset q^{-1}(C) \qand B' \subset q^{-1}(D),
\end{equation}
where $(C,D)$ is a Sturmian pair in $M$ such that $m_M(CD) = m_G(A' B')$. In particular, we have 
\[
m_G(A) = m_G(A') \leq m_M(C) \qand m_G(B) = m_G(B') \leq m_M(D). 
\]
However, since both $(A,B)$ and $(C,D)$ are sub-critical, 
we also have
\[
m_G(A) + m_G(B) = m_G(AB) \geq m_G(A'B') = m_M(CD) = m_M(C) + m_M(D),
\]
which now forces $m_G(A) = m_M(C)$ and $m_G(B) = m_M(D)$, and so in particular $m_G(AB) = m_M(CD)$. \\

We stress that it is not clear at this point whether $(A,B)$ reduces to $(C,D)$, that is to say, we do not
yet know that the inclusions in \eqref{defCD} also hold when $A'$ and $B'$ are replaced with $A$ and $B$ respectively. Even 
if this were true, we would still need an argument to show that this implies that the pair $(A_o, B_o)$ reduces to $(C,D)$, which is what Theorem \ref{main} asserts.
In order to fill in these gaps, we shall utilize the notions of \emph{stability} and \emph{regularity} of pairs of Borel sets. 

\subsubsection{Stability and regularity}

\begin{definition}[Regular and stable pairs]
\label{defreg}
We say that a closed set $A \subset G$ is \textbf{regular} if it is Jordan measurable with respect to $m_G$ and 
equal to the closure of its interior, and we say that a pair $(A,B)$ is \emph{regular} if both $A$ and $B$ are regular sets. We say that $(A,B)$ is \textbf{stable} if the 
inclusions 
\[
aB \subset AB \qand Ab \subset AB
\]
imply that $a \in A$ and $b \in B$.
\end{definition}

\begin{remark}
We leave it to the reader to verify that Sturmian pairs are always regular and stable, as is the the pair $(q^{-1}(C),
q^{-1}(D))$ in \eqref{defCD} above. 
\end{remark}

If $A_1$ and $A_2$ are Borel sets, we write $A_1 \sim A_2$ if 
$m_G(A_1 \Delta A_2) = 0$, where $\Delta$ denotes the symmetric difference of sets.

\begin{lemma}
\label{rigidity}
Suppose that $(A_1,B_1)$ and $(A_2,B_2)$ are sub-critical pairs in $G$ such that
\[
A_1\sim A_2 \qand B_1 \sim B_2.
\]
If $(A_2, B_2)$ is a regular and stable pair, then $A_1 \subset A_2$ and $B_1 \subset B_2$.
\end{lemma}

We recall from above that $A = p^{-1}(A_o)$ and $B = p^{-1}(B_o)$ and 
\[
p^{-1}(A_o) \sim q^{-1}(C) \qand p^{-1}(B_o) \sim q^{-1}(D),
\]
and both 
\[
(p^{-1}(A_o),p^{-1}(B_o)) \qand (q^{-1}(C),q^{-1}(D))
\]
are sub-critical pairs in $G$. The latter pair is in addition both regular and stable. Lemma \ref{rigidity}, 
applied to  
\[
A_1 =  p^{-1}(A_o) \qand A_2 = q^{-1}(C) \qand B_1 = p^{-1}(B_o)  \qand B_2 = q^{-1}(D),
\]
now tells us that
\begin{equation}
\label{defABCD}
 p^{-1}(A_o) \subset q^{-1}(C) \qand p^{-1}(B_o)  \subset q^{-1}(D).
\end{equation}
We claim that these inclusions force $\ker p \subset \ker q$, which will finish the proof of Theorem 
\ref{main}. Indeed, if  $\ker p \subset \ker q$, then the map $\pi : K \ra M$ given by 
$\pi(k) = q(p^{-1}(k))$ is a well-defined homomorphism, and by \eqref{defABCD}, we have
\[
 A_o \subset \pi^{-1}(C) \qand B_o  \subset \pi^{-1}(D).
\]
Furthermore, since $m_K(A_o B_o) = m_G(AB) = m_M(CD)$, the pair $(A_o,B_o)$ reduces to 
$(C,D)$, and 
\[
m_K(A_o) = m_G(A) = m_M(C) \qand m_K(B_o) = m_G(B) = m_M(D).
\]
\subsubsection{Proving $\ker p \subset \ker q$}

Let us return to \eqref{defABCD}. We first note that
\[
p^{-1}(A_o) \subset q^{-1}(C)g, \quad \textrm{for all $g \in \ker p$},
\]
whence
\[
p^{-1}(A_o) \subset q^{-1}(C) \cap q^{-1}(C)g = q^{-1}(C \cap Cq(g)),
\]
for all $g \in \ker p$, and thus 
\[
m_K(A_o) \leq m_M(C \cap Cq(g)) \leq m_M(C) = m_K(A_o), 
\]
which implies that $m_M(C) = m_M(C \cap C q(g))$ for all $g \in \ker p$. Hence $q(\ker p)$ is contained in the 
\emph{right (essential) stabilizer $Q$ of $C$}, where
\begin{equation}
\label{definQ}
Q = \big\{ m \in M \, : \, m_M(C \cap Cm) = m_M(C) \big\} < M.
\end{equation}
If denote by $r$ the right regular representation of $M$ on $L^2(M)$, then we see that $Q$ is the actual
stabilizer of the indicator function $\chi_C$ in $L^2(M)$. Since $r$ acts (norm-)continuously
on $L^2(M)$, we conclude that $Q$ is a \emph{closed} subgroup of $M$. \\

If $M = \bT$ and $C = aI$ for some $a \in \bT$ and \emph{proper} closed interval $I \subset \bT$, then $Q$ is clearly trivial, and 
thus $q(\ker p)$ is trivial as well.  If $M = \bT \rtimes \{-1,1\}$ and $C = a(I \rtimes \{-1,1\})$ for some $a \in \bT \rtimes \{-1,1\}$ and \emph{proper} closed interval $I \subset \bT$, then a tedious, but straightforward calculation shows that $Q = \{0\} \rtimes \{-1,1\}$, which is \emph{not} a normal subgroup of $M$. However, since the image $q(\ker p)$ is always a normal subgroup of $M$, which must be contained in $Q$, we conclude that the subgroup $q(\ker p)$ is trivial. 

\subsection{Proof of Lemma \ref{rigidity}}

We assume that $(A_1,B_1)$ and $(A_2,B_2)$ are sub-critical pairs in $G$ and $A_1 \sim A_2$ and $B_1 \sim B_2$, and $(A_2,B_2)$ is a regular and stable pair. We define the sets
\[
A_o = A_1 \cap A_2 \qand B_o = B_1 \cap B_2
\]
and note that
\[
m_G(A_o) = m_G(A_1) = m_G(A_2) 
\qand  
m_G(B_o) = m_G(B_1) = m_G(B_2),
\]
and
\begin{eqnarray}
\label{seqofeqs}
m_G(A_oB_o) 
&\leq &
m_G(A_1B_1) \nonumber \\
&= &
m_G(A_1) + m_G(B_1) \nonumber \\
&= &
m_G(A_o) + m_G(B_o) < 1. 
\end{eqnarray}
If the first inequality were strict, then $(A_o,B_o)$ is critical, which by 
the first assertion of Theorem \ref{knes} implies that there exists an open 
\emph{normal} subgroup $U \triangleleft G$ such that $A_o B_o = A_o B_o U$. We
note that $A_oU$ and $B_o U$ are clopen (closed and open) sets, and thus 
\begin{equation}
\label{V1V2}
V_1 = A_2^o \setminus A_o U \qand V_2 = B_2^o \setminus B_o U.
\end{equation}
are open. Since $A_2$ and $B_2$ are Jordan measurable and $A_o \sim A_2$
and $B_o \sim B_2$, we have $A_o \sim A^{o}_2$ and $B_o \sim B^{o}_2$. Hence, $V_1$ and $V_2$ are 
open $m_G$-null sets, and thus empty. Furthermore, since $A_2$ and $B_2$ are regular sets and $A_o U$
and $B_o U$ are clopen, 
we conclude by \eqref{V1V2} that
\[
A_2 = \overline{A_2^{o}} \subset A_o U \qand B_2 = \overline{B_2^{o}} \subset B_o U,
\]
and, since $U$ is normal, 
\[
A_2 B_2 \subset A_o U B_o = A_o B_o \subset A_2 B_2,
\]
and thus
\[
m_G(A_o B_o) = m_G(A_2 B_2).
\]
Since $(A_2, B_2)$ is sub-critical, we have
\begin{equation}
\label{chain2}
m_G(A_o B_o) = m_G(A_2 B_2) = m_G(A_2) + m_G(B_2) = m_G(A_o) + m_G(B_o) < 1,
\end{equation}
which contradicts our assumption that $(A_o,B_o)$ is critical. Going back to the chain of identities in 
\eqref{seqofeqs} and \eqref{chain2}, we see that we can henceforth assume that 
\[
m_G(A_oB_o) = m_G(A_1B_1) = m_G(A_2B_2).
\]
We wish to prove that $A_1 \subset A_2$ and $B_1 \subset B_2$. Assume for the sake of contradiction 
that there exists an element $x \in A_1 \setminus A_2$. In particular, $x \notin A_2$ and $A_o \cup \{x\} \subset A_1$. Then, 
since $A_o \sim A_1 \sim A_2$ and $B_o \sim B_1 \sim B_2$ and $A_o B_o \sim A_2 B_2$,
we have
\begin{eqnarray*}
m_G(A_1) + m_G(B_1) 
&= &
m_G(A_1B_1) \\
&\geq &
m_G((A_o \cup \{x\})B_o) \\
&= &
m_G(A_oB_o \cup (xB_o \setminus A_oB_o)) \\
&= &
m_G(A_o B_o) + m_G(xB_o \setminus A_oB_o) \\
&= &
m_G(A_o) + m_G(B_o) + m_G(xB_2 \setminus A_2B_2)) \\
&= &
m_G(A_1) + m_G(B_1) + m_G(xB_2 \setminus A_2B_2),
\end{eqnarray*}
which forces the set $xB_2 \setminus A_2 B_2$ to be a $m_G$-null set. 
Since $A_2$ and $B_2$ are closed sets, so is $A_2 B_2$. 
In particular, the set $x B_2^o \setminus A_2 B_2$ is open 
and $m_G$-null, which forces it to be empty. Since $B_2$ is regular and $A_2 B_2$ is closed, 
we have 
\[
x\overline{B_2^o} = xB_2 \subset A_2 B_2.
\]
By assumption, the pair $(A_2,B_2)$ is stable, and thus $x \in A_2$, 
which is a contradiction. This shows that $A_1 \subset A_2$. The proof
that $B_1 \subset B_2$ works the same.

\section{Proof of Proposition \ref{relativization}}
\label{proofrel}

Let $G$ be a compact and second countable group with identity component $N$. Since $G$ is second countable,
there exists a decreasing sequence $(U_n)$ of \emph{open} normal subgroups of $G$ such that $N$ equals their 
intersection. The following two propositions immediately imply Proposition \ref{relativization}.

\begin{proposition}
\label{criticalU}
Suppose that $(A,B)$ is spread-out and sub-critical in $G$. Then $(A,B)$ is sub-critical with 
respect to $U_n$ for every $n$.
\end{proposition}

\begin{proposition}
\label{criticalN}
If $(A,B)$ is sub-critical with respect to $U_n$ for every $n$, then $(A,B)$ is sub-critical with respect
to $N$.
\end{proposition}

\subsection{Proof of Proposition \ref{criticalU}}
Let us fix an open normal subgroup $U$ of $G$ throughout this subsection, and 
suppose that $(A,B)$ is spread-out and sub-critical in $G$. We note that the Haar 
probability measure $m_{U}$ on the clopen subgroup $U$, viewed as a Borel probability 
measure on $G$, is given by
\begin{equation}
\label{defmUn}
m_{U}(D) = \frac{m_G(D \cap U)}{m_G(U)}, \quad \textrm{for Borel sets $D \subset G$}.
\end{equation}
We wish to prove that $(A,B)$ is sub-critical with respect to $U$. Given a Borel set $D \subset G$, 
we define 
\[
D_x = D \cap xU, \quad \textrm{for $x \in G$}.
\] 
We note that $D_x$ only depends on the right $U$-coset of $x$,  so we may just as well view $x$ as an element in $G/U$. Furthermore, 
\[
A_x = x(x^{-1}A \cap U) \qand B_y = (By^{-1} \cap U)y
\]
and
\[
A_x B_y = x(x^{-1}A \cap U)(By^{-1} \cap U)y \subset AB \cap Uxy = (AB)_{xy}
\]
for all $x, y \in G$. Using \eqref{defmUn}, the left and right invariance of the Haar probability measure $m_G$ and 
the identity $m_G(U) = \frac{1}{|G/U|}$, 
we see that the conditions for sub-criticality of $(A,B)$ with respect to $U$ (see Definition \ref{defsubcrit}) can be equivalently rewritten as:
\begin{equation}
\label{show1}
m_G(A_x)= \frac{m_G(A)}{|G/U|} \qand m_G(B_y) = \frac{m_G(B)}{|G/U|}
\end{equation}
and 
\begin{equation}
\label{show2}
m_G((AB)_{xy}) = m_G(A_x B_y) = m_G(A_x) + m_G(B_y) < m_G(U),
\end{equation}
for all $x,y \in G/U$. We stress that the last inequality in \eqref{show2} follows from \eqref{show1} and
our assumption that $(A,B)$ is sub-critical. \\

To prove these identities, the following technical lemma will be very useful.

\begin{lemma}
\label{mxy}
For all $x, y \in G/U$, the sets $A_x$ and $B_y$ are non-empty, and 
\[
m_G((AB)_{xy}) \geq m_G(A_x) + m_G(B_y).
\]
\end{lemma}

Before we prove this lemma, we show how to deduce \eqref{show1} and \eqref{show2} from it. We first note
that
\[
AB = \bigsqcup_{z \in G/U} (AB)_z = \bigsqcup_{z \in G/U} \Big( \bigcup_{xy = z} A_x B_y \Big).
\]
Pick $x_o, y_o \in G/U$ such that
\[
m_G(A_{x_o}) = \max_{x \in G/U} m_G(A_x) \qand m_G(B_{y_o}) = \max_{y \in G/U} m_G(B_y),
\]
and note that
\begin{equation}
\label{Axo}
m_G(A) \leq |G/U| m_G(A_{x_o}) \qand m_G(B) \leq |G/U| m_G(B_{y_o}).
\end{equation}
By Lemma \ref{mxy}, we have
\[
m_G((AB)_z) \geq m_G(A_{zy_o^{-1}}) + m_G(B_{y_o})
\qand
m_G((AB)_z) \geq m_G(A_{x_o}) + m_G(B_{x_o^{-1}z}),
\]
for all $z \in G/U$. Hence,
\[
m_G(AB) \geq \sum_{z \in G/U} \big( m_G(A_{zy_o^{-1}}) + m_G(B_{y_o})\big) = m_G(A) + |G/U| m_G(B_{y_o})
\]
and
\[
m_G(AB) \geq \sum_{z \in G/U} \big( m_G(A_{x_o}) + m_G(B_{x_o^{-1} z})\big) = |G/U| m_G(A_{x_o}) + m_G(B).
\]
By \eqref{Axo} and sub-criticality of $(A,B)$, we see that
\[
m_G(A_{x_o}) = \frac{m_G(A)}{|G/U|} \qand m_G(B_{y_o}) = \frac{m_G(B)}{|G/U|}. 
\]
This implies that
\[
m_G(A) = \sum_{x \in G/U} m_G(A_x) \leq |G/U| m_G(A_{x_o}) = m_G(A),
\]
and similarly for $B$. Hence,
\[
m_G(A_{x}) = \frac{m_G(A)}{|G/U|} \qand m_G(B_{y}) = \frac{m_G(B)}{|G/U|}, \quad \textrm{for all $x, y \in G/U$}.
\]
This proves \eqref{show1}. We can now replace $x_o$ and $y_o$ with arbitrary $x$ and $y$ in the inequalities above. Using Lemma \ref{mxy}, we conclude that for every fixed $y$, we have
\begin{eqnarray*}
m_G(AB) 
&\geq& \sum_{z \in G/U} m_G((AB)_z) 
\geq \sum_{z \in G/U} m_G(A_{zy^{-1}}) + m_G(B_y) \\
&=& m_G(A) + m_G(B).
\end{eqnarray*}
Since $(A,B)$ is sub-critical, these inequalities are in fact equalities, and we now see that
\[
m_G((AB)_z) = m_G(A_{zy^{-1}}B_y) = m_G(A_{zy^{-1}}) + m_G(B_y), \quad \textrm{for all $y,z \in G/U$}.
\]
This shows \eqref{show2}, and thus the proof of Proposition \ref{criticalU} is finished. 

\subsection{Proof of Lemma \ref{mxy}}

First note that if $A$ is spread-out in $G$, then $A_x$ is non-empty for every $x \in G/U$. Indeed, if it were empty
for some $x$, then $xU \notin q(A)$ where $q$ denotes the canonical quotient map onto $G/U$, and thus $A$ does not project onto the finite group $G/U$, which contradicts our assumption that $A$ is spread-out. \\

If $X$ is a subset of $G$, we denote by $X^c$ the complement of $X$ in $G$. Suppose that $(A,B)$ is sub-critical in $G$ and define $C = (AB)^c$. Then $A^{-1}C \subset B^c$, and thus
\[
m_G(A^{-1}C) \leq 1 - m_G(B) = m_G(A) + m_G(C) < 1,
\]
since $B$ is not a null set. We note that if the first inequality is strict, then the pair $(A^{-1},C)$ is critical in $G$, and thus
$A$ cannot be spread-out by the first assertion in Theorem
\ref{knes}. We conclude that $(A^{-1},C)$ is sub-critical in $G$. Furthermore, 
\begin{equation}
\label{Bxz}
(A^{-1})_x C_z \subset (A^{-1}C)_{xz} \subset (B^c)_{xz}, \quad \textrm{for all $x,z \in G/U$}.
\end{equation}
\begin{lemma}
\label{cxy}
For all $x, z \in G/U$, we have
\[
m_G((A^{-1})_x) + m_G(C_z) \leq m_G((A^{-1})_x C_z).
\]
\end{lemma}
We claim that this lemma implies Lemma \ref{mxy}. Indeed, first note that for every Borel set $D \subset G$, we 
have
\begin{equation}
\label{Dxinv}
(D^{-1})_z = D_{z^{-1}}^{-1} \qand m_G((D^c)_z) = m_G(U) - m_G(D_z), \quad \textrm{for all $z \in G/U$}.
\end{equation}
In particular, by \eqref{Bxz} and the definition of $C$,
\[
m_G((A^{-1})_x C_z) \leq m_G(U) - m_G(B_{xz}) 
\qand
m_G(C_z) = m_G(U) - m_G((AB)_z),
\]
for all $x, z \in G/U$. Suppose that Lemma \ref{cxy} holds. Then, using these relations, we conclude that
\[
m_G((A^{-1})_x) + m_G(U) - m_G((AB)_z) \leq m_G((A^{-1})_x C_z) \leq m_G(U) - m_G(B_{xz}),
\]
which readily translates to 
\[
m_G(A_{x^{-1}}) + m_G(B_{xz}) \leq m_G((AB)_z), \quad \textrm{for all $x, z \in G/U$},
\] 
where we have used the relation $(A^{-1})_x = A_{x^{-1}}^{-1}$ from \eqref{Dxinv}, and the fact that 
$m_G$ is inversion-invariant (since $G$ is compact, and thus unimodular). The proof of Lemma
\ref{mxy} is now complete. 

\subsubsection{Proof of Lemma \ref{cxy}}
Suppose that $(S,T)$ is a sub-critical pair in $G$ and suppose that
\begin{equation}
\label{assume}
m_G(S_x T_y) < m_G(S_x) + m_G(T_y), \quad \textrm{for \emph{some} $x, y \in G/U$}.
\end{equation}
This translates to the bound
\[
m_U((x^{-1}S \cap U)(Ty^{-1} \cap U)) < m_U(x^{-1}S \cap U) + m_U(Ty^{-1} \cap U),
\]
and thus we see that the pair $((x^{-1}S \cap U), (Ty^{-1} \cap U))$ is critical in $U$. There
are now two cases to consider. \\

\textbf{Case I:} Suppose that 
\[
m_U((x^{-1}S \cap U)(Ty^{-1} \cap U)) < 1.
\]
Then $((x^{-1}S \cap U), (Ty^{-1} \cap U))$ reduces to a pair of \emph{proper} subsets of a finite quotient
group $U/Q$ for some open proper normal subgroup $Q$ of $U$. In particular, 
\[
(x^{-1}S \cap U)Q  \neq U.
\]
We can also view $Q$ as an open (but not necessarily normal) subgroup of $G$. However, since $G/Q$ is 
finite, we see that $R = \bigcap_{g \in G/Q} gQg^{-1}$ is an open \emph{normal} subgroup of $G$, 
which by construction is contained in $Q$ and thus U$.$ In particular, 
\[
(x^{-1}S \cap U)R \neq U,
\]
whence
\[
x^{-1}SR = (x^{-1}S \cap U)R \cup (x^{-1}S \cap U^c)R \subset (x^{-1}S \cap U)R \cup U^c \neq G.
\]
We see that $S$ does not project onto the finite quotient group $G/R$, and thus $S$ is \emph{not} spread-out. \\

\textbf{Case II:} Suppose that
\[
1 = m_U((x^{-1}S \cap U)(Ty^{-1} \cap U)) < m_U(x^{-1}S \cap U) + m_U(Ty^{-1} \cap U).
\]
Then it is not hard to see that the product set equals $U$, and thus $S_x T_y = xyU \subset ST$. \\

Hence we have the following alternative: If $(S,T)$ is sub-critical and \eqref{assume} holds for some $x, y \in G/U$, then either
\begin{itemize}
\item $S$ is \emph{not} spread-out, or
\item $ST$ contains a coset of $U$.
\end{itemize}
Let us now apply this observation to the pair $(S,T) = (A^{-1},C)$ in Lemma \ref{cxy} above. Since $A$ is assumed to be
spread-out, so is $A^{-1}$, and thus the assertion in Lemma \ref{cxy} can only fail if $ST$ contains a coset of $U$, i.e. if $A^{-1}C \supset zU$ for some $z \in G/U$. However, recall that $A^{-1}C \subset B^c$. We conclude that $B \cap zU = \emptyset$, and thus $B$ is 
not spread-out, contrary to our assumption.

\subsection{Proof of Proposition \ref{criticalN}}
Let $(U_n)$ be a decreasing sequence of open and normal subgroups of $G$ with intersection $N$. We recall
that the Haar probability measures on $U_n$ can be viewed as Borel probability measures on $G$ via
\[
m_{U_n}(C) = \frac{m_G(C \cap U_n)}{m_G(U_n)}, \quad \textrm{for Borel sets $C \subset G$}.
\]
By uniqueness of Haar probability measures on compact groups, we observe that $m_{U_n} \ra m_N$ in the
weak*-topology on the space of Borel probability measures on $G$. Suppose that $(A,B)$ is a sub-critical pair 
with respect to $U_n$ for every $n$, and set $\mu_n = m_{U_n}$. We leave it as an exercise to show that since 
$U_n$ is open, we can take the sets $X$ and $Y$ in Definition \ref{defsubcrit} to be $G$ and $G \times G$.  \\

We note that \eqref{1} can be rewritten as
\begin{equation}
\label{Asy}
m_G(A) = \int_{G} \chi_A(sy) \, d\mu_n(y)
\qand
m_G(B) = \int_{G} \chi_{B}(yt) \, d\mu_n(y), 
\end{equation}
for all $s, t \in G$, and if we combine \eqref{1} with \eqref{2} and \eqref{3}, we have
\begin{equation}
\label{AB}
1 > m_G(A) + m_G(B) = \int_G \chi_{AB}(syt) \, d\mu_n(y),
\end{equation}
for all $s, t \in G$. \\

We claim that \eqref{Asy} and \eqref{AB} still hold for a \emph{conull} set of $(s,t) \in G \times G$ if $\mu_n$ 
is replaced with $m_N$. This will follow from Lemma \ref{mun} below. After the proof of this lemma, we will 
show how this can be used to finish the proof of Proposition \ref{criticalN}.

\begin{lemma}
\label{mun}
Let $(\mu_n)$ be a sequence of Borel probability measures on $G$ which converges in the weak*-topology
to a Borel probability measure $\mu$ on $G$. Then, for every bounded real-valued Haar measurable function 
$f$ on $G$, there exist a subsequence $(n_j)$ and conull subsets $X \subset G$ and $Y \subset X \times X$ 
such that 
\[
\int_G f(sy) \, d\mu_{n_j}(y) \ra \int_G f(sy) \, d\mu(y) 
\qand
\int_G f(yt) \, d\mu_{n_j}(y) \ra \int_G f(yt) \, d\mu(y)
\]
for all $s, t \in X$, and
\[
\int_G f(syt) \, d\mu_{n_j}(y) \ra \int_G f(syt) \, d\mu(y), \quad \textrm{for all $(s,t) \in Y$}.
\]
\end{lemma}

\begin{proof}
Since $\mu_n \ra \mu$ in the weak*-topology, the lemma is trivial when $f$ is continuous, and then 
no passage to a sub-sequence is necessary. By a standard approximation argument, combined with
dominated convergence, the lemma also holds for every bounded Haar measurable function on $G$ with respect to the 
\emph{norm} topologies on $L^2(G)$ and $L^2(G \times G)$ respectively. Since every $L^2$-convergent
sequence admits an almost everywhere convergent sub-sequence, we are done.
\end{proof}

Applied to the sequence $\mu_n = m_{U_n}$ above and $\mu = m_N$ and the relations \eqref{Asy} and \eqref{AB}, we conclude from Lemma \ref{mun} that there are conull 
Borel sets $X \subset G$ and $Y \subset X \times X$ such that
\begin{equation}
\label{AN}
m_G(A) = m_N(s^{-1}A \cap N) \qand m_G(B) = m_N(Bt^{-1} \cap N)
\end{equation}
for all $s,t \in X$, and
\begin{equation}
\label{ABN}
1 > m_G(A) + m_G(B) = m_N(s^{-1}ABt^{-1} \cap N), \quad \textrm{for all $(s,t) \in Y$}.
\end{equation}
Since $N$ is connected and abelian, we know by the second
assertion in Theorem \ref{knes} that the pairs $((s^{-1}A \cap N),(Bt^{-1} \cap N))$ are 
never critical, and since $s^{-1}ABt^{-1} \supset (s^{-1} A \cap N)(Bt^{-1} \cap N)$ we have
\[
m_N(s^{-1}ABt^{-1} \cap N) \geq m_N((s^{-1} A \cap N)(Bt^{-1} \cap N)) \geq m_N(s^{-1}A \cap N) + m_N(Bt^{-1} \cap N).
\]
Upon combining \eqref{AN} and \eqref{ABN}, we now see that
\[
m_N(s^{-1}ABt^{-1} \cap N) = m_N((s^{-1} A \cap N)(Bt^{-1} \cap N)) 
\]
and
\[
m_N((s^{-1} A \cap N)(Bt^{-1} \cap N)) = m_N(s^{-1}A \cap N) + m_N(Bt^{-1} \cap N),
\]
for all $(s,t) \in Y$, which shows that $(A,B)$ is sub-critical with respect to $N$.

\section{Proof of Proposition \ref{relativemain}}
\label{proofrelmain}

Let $G$ be a compact and second countable group with an \emph{abelian} identity component $N$ and a closed 
subgroup $L$ such that $N \cap L = \{e\}$ and $NL = G$. Let $A, B \subset G$ be 
Borel sets and suppose that $(A,B)$ is sub-critical with respect to $N$. We recall that this means that there
are conull subsets $X \subset G$ and $Y \subset X \times X$ such that
\begin{equation}
\label{subAandB}
m_N(s^{-1}A \cap N) = m_G(A) \qand m_N(Bt^{-1} \cap N) = m_G(B),
\end{equation}
for all $s, t \in X$, and
\begin{equation}
\label{subAorB}
m_N((s^{-1} A \cap N)(Bt^{-1} \cap N)) = m_G(A) + m_G(B) < 1
\end{equation}
and 
\begin{equation}
\label{subAB}
m_N(s^{-1}ABt^{-1} \cap N) = m_N((s^{-1} A \cap N)(Bt^{-1} \cap N)),
\end{equation}
for all $(s,t) \in Y$.  Since $m_G$ is inversion-invariant, we may without loss of generality assume that $X$ and $Y$ 
are invariant under taking inverses, so in particular, the identities \eqref{subAandB}, \eqref{subAorB}
and \eqref{subAandB} above also hold with $s^{-1}$ replaced with $s$. We shall henceforth 
assume that these replacements have been made. 

\subsection{A basic reduction}

Fix $s,t \in G$. Since $N \cap L = \{e\}$ and $NL = G$, we can write $s = n_s l_s$ and $t = n_t l_t$
for unique elements $n_s, n_t \in N$ and $l_s, l_t \in L$. Hence, if $s,t \in X$ and $(s,t) \in Y$, then
\[
sA \cap N = n_s(l_s A \cap N) \qand Bt^{-1} \cap N = (B l_t^{-1} \cap N) n_t^{-1}
\]
and 
\[
sABt^{-1} \cap N= n_s(l_s AB l_t^{-1} \cap N)n_t^{-1}.
\] 
Since $m_N$ is left and right $N$-invariant, we conclude
that the relations \eqref{subAandB}, \eqref{subAorB} and \eqref{subAB} (with $s$ replaced with $s^{-1}$) only depend on $l_s$ and $l_t$.
In particular, the sets $X$ and $Y$ are respectively left $N$- and $N \times N$-invariant, so we conclude that there are conull sets $X_o \subset L$ and $Y_o \subset L \times L$
such that
\begin{equation}
\label{subAandBL}
m_N(sA \cap N) = m_G(A) \qand m_N(Bt^{-1} \cap N) = m_G(B),
\end{equation}
for all $s, t \in X_o$, and
\begin{equation}
\label{subAorBL}
m_N((sA \cap N)(Bt^{-1} \cap N)) = m_G(A) + m_G(B) < 1
\end{equation}
and 
\begin{equation}
\label{subABL}
m_N(sABt^{-1} \cap N) = m_N((s A \cap N)(Bt^{-1} \cap N)),
\end{equation}
for all $(s,t) \in Y_o$. The main difference from \eqref{subAandB}, \eqref{subAorB} and \eqref{subAB} is that
$s$ and $t$ are now elements of $L$.

\subsection{Combing the dual group of $N$}
\label{subseccomb}
Recall that $N$ is assumed to be a compact and \emph{connected} abelian group. Let $\widehat{N}$ denote the 
dual group of $N$, i.e. the group of all continuous homomorphisms from $N$ to $\bT$ with pointwise
addition (which we write multiplicatively).  It is a classical fact that $\widehat{N}$ is a countable and \emph{torsion-free} group. In particular, the map $\xi \mapsto \check{\xi}$ on $\widehat{N}$ given by
\[
\check{\xi}(n) = \xi(n)^{-1}, \quad \textrm{for $n \in N$},
\]
has only one fixed point, namely the trivial homomorphism, here denoted by $1$. Hence, we can inductively construct a set $S \subset \widehat{N} \setminus \{1\}$ with the property that
\begin{equation}
\label{defS}
\widehat{N} \setminus \{1\} = S \cup \check{S} \qand S \cap \check{S} = \emptyset.
\end{equation}
Suppose that we have fixed such a set $S$ once and for all. \\

Note that \eqref{subAandBL} and \eqref{subAorBL} above show that
\[
(sA \cap N, Bt^{-1} \cap N), \quad \textrm{for $y = (s, t) \in Y_o$}
\]
are all sub-critical pairs in $N$. Since $N$ is abelian, Theorem \ref{knes} asserts
they all reduce to Sturmian pairs in $\bT$. Recall that this means that for every 
$y = (s,t) \in Y_o$, we can find \\
\begin{itemize}
\item a continuous homomorphism $\xi_y : N \ra \bT$. \\
\item closed and symmetric intervals $I_y, J_y \subset \bT$ such that 
\begin{equation}
\label{IyJy}
m_N((sA \cap N)(Bt^{-1} \cap N)) = m_{\bT}(I_y J_y).
\end{equation}
\item $a(y), b(y) \in \bT$ such that
\[
sA \cap N \subset \xi_y^{-1}(I_y a(y)) \qand Bt^{-1} \cap N \subset \xi_y^{-1}(J_y b(y)).
\]
In particular, by \eqref{subAandBL}, 
\begin{equation}
\label{AlessI}
m_G(A) \leq m_{\bT}(I_y) \qand m_G(B) \leq m_{\bT}(J_y).
\end{equation}
\end{itemize}

\begin{remark} 
Since $I_y$ and $J_y$ are symmetric, we have
\[
\xi_y^{-1}(I_ya(y)) = \check{\xi}_y^{-1}(I_ya(y)^{-1}) 
\qand
\xi_y^{-1}(I_yb(y)) = \check{\xi}_y^{-1}(I_yb(y)^{-1}).
\quad \textrm{for all $y$},
\]
Hence we can, possibly upon changing $a(y)$ and $b(y)$ to $a(y)^{-1}$ and $b(y)^{-1}$ whenever necessary, 
assume that $\xi_y \in S$ for all $y \in Y_o$. We shall make this assumption throughout the rest of the paper. 
\end{remark}

Furthermore, since $(I_y, J_y)$ is clearly sub-critical in $\bT$ for every $y \in Y_o$, we have by \eqref{subAorBL}
and \eqref{IyJy},
\[
m_G(A) + m_G(B) = m_N((sA \cap N)(Bt^{-1} \cap N)) = m_{\bT}(I_y J_y) = m_{\bT}(I_y) + m_{\bT}(J_y) < 1.
\]
By \eqref{AlessI}, we conclude that $m_G(A) = m_{\bT}(I_y)$ and $m_G(B) = m_{\bT}(J_y)$. Since we have 
assumed that $I_y$ and $J_y$ are both closed and \emph{symmetric} intervals, they are \emph{uniquely} determined by their Haar 
measures. In particular, we conclude that $I_y$ and $J_y$ are independent of $y$, and we shall henceforth
simply denote them by $I$ and $J$. \\

\emph{To summarize:} Let $I, J \subset \bT$ denote the \emph{unique} closed and symmetric intervals of $\bT$
of Haar measures $m_G(A)$ and $m_G(B)$ respectively. Fix a set $S \subset \widehat{N} \setminus \{1\}$ as in
\eqref{defS}. Then, for every $y = (s,t) \in Y_o$, there exist $\xi_y \in S$ and $a(y), b(y) \in \bT$ such that 
\begin{equation}
\label{xiy}
sA \cap N \subset \xi_y^{-1}(I a(y)) \qand Bt^{-1} \cap N \subset \xi_y^{-1}(J b(y)).
\end{equation}
Furthermore, we have
\begin{equation}
\label{AsN}
m_N(sA \cap N) = m_{\bT}(I) \qand m_N(Bt^{-1} \cap N) = m_{\bT}(I).
\end{equation}

\subsection{Getting rid of dependencies}
\label{getrid}

If $E, F \subset N$ are Borel sets, we write
$E \sim F$ if $m_N(E \Delta F) = 0$, where $\Delta$ denotes the symmetric difference of sets.  \\

By \eqref{xiy} and \eqref{AsN}, we see that whenever $y = (s,t)$ and $y' = (s,t')$ are elements having the same
\emph{first} coordinate, and which both belong to $Y_o$, then
\[
\xi_{y}^{-1}(Ia(y)) \sim \xi_{y'}^{-1}(Ia(y')).
\]
It is not hard to see (using the fact that $I$ is regular and has trivial stabilizer in $\bT$), that this forces either
\[
\xi_y = \xi_{y'} \qand a(y) = a(y') 
\quad
\textbf{or}
\quad 
\xi_y = \check{\xi}_{y'} \qand a(y) = a(y')^{-1}.
\]
However, since $\xi_y, \xi_{y'} \in S$ for all $y, y' \in Y_o$, and $S \cap \check{S} = \emptyset$, the second possibility cannot occur. Similarly, 
if $z = (u,v)$ and $z' = (u',v)$ are elements having the same
\emph{second} coordinate, and which both belong to $Y_o$, then
\[
\xi_{z}^{-1}(Jb(z)) \sim \xi_{z'}^{-1}(Jb(z')),
\]
which forces either 
\[
\xi_z = \xi_{z'} \qand b(z) = b(z') 
\quad
\textbf{or}
\quad 
\xi_z = \check{\xi}_{z'} \qand b(z) = b(z')^{-1}.
\]
Since $\xi_{z}, \xi_{z'} \in S$ for all $z,z ' \in Y_o$, the second possibility does not occur, and we conclude that 
\begin{equation}
\label{cool}
\xi_y = \xi_{y'} \qand \xi_z = \xi_{z'} \qand a(y) = a(y') \qand b(z) = b(z').
\end{equation}
We claim that these identities imply that there exist $\xi \in S$ such that $\xi_y = \xi$ for a.e. $y$, and
functions $\alpha, \beta : X_o \ra \bT$ (possibly upon shrinking $X_o$ to a conull subset thereof) such that
\begin{equation}
\label{abalpha}
a(y) = \alpha(s) \qand b(y) = \beta(t), \quad \textrm{for a.e. $y = (s,t)$}.
\end{equation}

\noindent \textbf{Proof of claim:} First note that the continuous map $q : L^4 \ra L^2$ defined by
\[
q((s,t),(u,v)) = (s,v), \quad \textrm{for $(s,t), (u,v) \in L^2$},
\]
maps the Haar measure on $L^4$ to the Haar measure on $L^2$. In particular, since $Y_o \subset L^2$
is a conull Borel set, the set
\[
F = q^{-1}(Y_o) \cap (Y_o \times Y_o) = \big\{ ((s,t),(u,v)) \in Y_o \times Y_o \, : \, (s,v) \in Y_o \big\} 
\]
is a conull Borel set of $Y_o \times Y_o$, so by Fubini's Theorem, there exists $(s,t) \in Y_o$ such that
the section 
\[
F_{(s,t)} = \big\{ (u,v) \in Y_o \, : \, ((s,t),(u,v)) \in F \big\} 
\]
is conull. Set $\xi = \xi_{(s,t)}$, and pick $(u,v) \in F_{(s,t)}$. By construction we have
$(s,t), (s,v), (u,v) \in Y_o$, so by \eqref{cool}, we must have
\[
\xi = \xi_{(s,t)} = \xi_{(s,v)} = \xi_{(u,v)}.
\]
In other words, $\xi_{(u,v)} = \xi$ for almost every $(u,v) \in Y_o$, which proves the first assertion. 
To prove \eqref{abalpha}, we argue as follows. By Theorem A.9 in \cite{Zi}, upon possibly replacing
$X_o$ with a conull Borel subset thereof, we can find Borel maps $q_1, q_2 : X_o \ra L$ such that
\[
(u,q_1(u)) \in Y_o \qand (q_2(u),u) \in Y_o, \quad \textrm{for all $u \in X_o$}.
\]
If we now define the (not a priori Borel measurable) maps $\alpha, \beta : X_o \ra \bT$ by
\[
\alpha(u) = a(u,q_1(u)) \qand \beta(u) = b(q_2(u),u), \quad \textrm{for $u,v \in X_o$},
\]
then by \eqref{cool}, we see that for all $(u,v) \in Y_o \cap (X_o \times X_o)$,
\[
a(u,v) = a(u,q_2(u)) = \alpha(u) \qand b(u,v) = b(q_2(v),v) = \beta(v),
\]
and thus we have proved \eqref{abalpha}. 

\begin{remark}
Since the Borel measurability of $q_1$ and $q_2$ is irrelevant at this point, some readers might wish to 
\emph{not} invoke Theorem A.9 in \cite{Zi} to prove the existence of $q_1$ and $q_2$. Instead, one can first 
extract a common conull Borel subset of $X_o$ (which we henceforth identify with $X_o$) of the projections 
of $Y_o$ onto each coordinate axis, and then use the 
Axiom of Choice to produce right inverses $q_1$ and $q_2$ of the coordinate projections restricted to $Y_o \cap (X_o \times X_o)$. 
\end{remark}

\emph{To summarize:} There are conull Borel subsets $X_o \subset L$ and $Y_o \subset X_o \times X_o$ 
(which may be different from the sets $X_o$ and $Y_o$ in the beginning of this subsection), maps 
$\alpha, \beta : X_o \ra \bT$ (with no obvious regularity whatsoever) and $\xi \in S$ such that
\begin{equation}
\label{summary}
sA \cap N \subset \xi^{-1}(I \alpha(s)) 
\qand
Bt^{-1} \cap N \subset \xi^{-1}(J \beta(t)), 
\end{equation}
for all $(s,t) \in Y_o$, where $I$ and $J$ denote the \emph{unique} closed and symmetric intervals in $\bT$ with 
Haar measures equal to $m_G(A)$ and $m(B)$ respectively. Furthermore, we also have
\[
m_N(sA \cap N) = m_{\bT}(I) \qand m_N(Bt^{-1} \cap N) = m_{\bT}(J).
\]
It now follows from \eqref{subABL}, \eqref{IyJy} and \eqref{xiy}, combined with the observation
that $\xi_y = \xi$ almost everywhere and \eqref{abalpha}, that 
\begin{equation}
\label{summaryAB}
sABt^{-1} \cap N \sim (sA \cap N)(Bt^{-1} \cap N) \sim \xi^{-1}(IJ \alpha(s) \beta(t)) 
\end{equation}
for all $(s,t) \in Y_o$. Note that we may without loss of generality (upon further restrictions) assume that $X_o$ is a symmetric subset of 
$G$ (since $m_G$ is inversion-invariant), and that the projections of $Y_o$ onto each $L$-coordinate coincide with $X_o$. We shall henceforth assume
that these assumptions are satisfied. \\

\emph{Technical interludes:} In what follows, we shall prove that $\xi \in S$ and the maps $\alpha, \beta : X_o \ra \bT$ can be used to construct a \emph{continuous}
homomorpism $\pi$ from $G$ into $\bT \rtimes \{-1,1\}$ so that the sets $A$ and $B$, modulo null sets, are contained in pre-images of a Sturmian pair in $\bT \rtimes \{-1,1\}$ under this homomorphism. Note however that we have yet not established any regularity, not even measurability, for the maps $\alpha$ and $\beta$. The technical tools for this will be outlined in the next two subsections. 

\subsection{Interlude I: Borel measurability of $\alpha$ and $\beta$}

We note that if $M$ is a compact group with Haar probability measure $m_M$ and $D \subset M$ is
a Borel set, then
\[
\Stab_M(D) = \big\{ m \in M \, : \, m_M(D m \cap D) = m_M(D) \big\}
\]
is a closed subgroup of $M$ (see the discussion after \eqref{definQ}). We see that if $M = \bT$ and $I$ is a proper closed interval of $\bT$,
then $\Stab_{\bT}(I)$ is the trivial subgroup. 

\begin{lemma}
\label{step4}
Suppose that 
\begin{itemize}
\item $G$ and $M$ are compact and second countable groups, and $L, N < G$ are closed subgroups. 
\item $\xi : N \ra M$ is a surjective continuous homomorphism.
\item $C \subset G$ and $I \subset M$ are Borel sets, and $\Stab_M(I)$ is trivial. 
\item $X \subset L$ is conull and there exists a map $\gamma : X \ra M$ such that  
\[
sC \cap N \sim \xi^{-1}(I \gamma(s)), \quad \textrm{for all $s \in X$}.
\]
\end{itemize} 
Then $\gamma$ is Borel measurable.
\end{lemma}

We shall use this lemma as follows. Since $I$ and $J$ are symmetric, we have by \eqref{AsN} and 
\eqref{summary} that
\[
sA \cap N \sim \xi^{-1}(I\alpha(s))
\qand
tB^{-1} \cap N \sim \xi^{-1}(J\beta(t)^{-1}),
\]
for all $s, t \in X_o$.  Applied to $G, N$ and $L$ as in the previous 
subsections, and $X = X_o \subset L$ and $M = \bT$ and 
\[
C = I \qand \gamma(s) = \alpha(s) \qor C = J \qand \gamma(s) = \beta(s)^{-1},
\]
the lemma above implies that $\alpha$ and $\beta$ are in fact Borel measurable as maps from $X_o \ra \bT$.

\begin{proof}[Proof of Lemma \ref{step4}]
Fix a countable basis $(U_n)$ for the topology on $M$ and note that, by assumption, we have
\begin{equation}
\label{eqweird}
sC \cap N \cap \xi^{-1}(U_n) \sim \xi^{-1}(I \gamma(s) \cap U_n)
\end{equation}
for all $n$ and for all $s$ in $X$. Define the maps
\[
\Psi : M \ra [0,1]^{\bN} \qand \Phi : L \ra [0,1]^{\bN} 
\]
by
\[
\Psi(t)_n = m_M(I \, t \cap U_n)
\qand
\Phi(s)_n = m_N\big(sC \cap N \cap \xi^{-1}(U_n)\big)
\]
for $n \geq 1$ and $t \in M$ and $s \in L$. We claim that both $\Psi$ and $\Phi$ are Borel 
measurable. To prove this, it suffices to show that $\Psi(\cdot)_n$ and 
$\Phi(\cdot)_n$ are Borel measurable for every $n$. Note that $M$ and $L$ act jointly 
continuously on $M$ and $N$ respectively, and thus they also act jointly continuously 
on the space of Borel probability measures on $M$ and $N$ respectively, endowed with
the weak*-topology. Hence, 
\[
t \mapsto \int_{U_n} f_1(mt) \, dm_M(m) \qand s \mapsto \int_{\xi^{-1}(U_n)} f_2(s^{-1}x) \, dm_N
\]
are continuous functions on $M$ and $L$ respectively, for every fixed pair $(f_1,f_2)$ of continuous 
function on $M$ and $N$ respectively. If we instead plug in $f_1 = \chi_I$ and $f_2 = \chi_C$, then 
these functions coincide with $\Psi(\cdot)_n$ and $\Phi(\cdot)_n$ respectively. Since both $f_1$ and
$f_2$ are pointwise limits of sequences of continuous functions, it follows from monotone convergence that both 
$\Psi(\cdot)_n$ and $\Phi(\cdot)_n$ are pointwise limits of sequences of continuous functions, and thus Borel 
measurable. \\

We further claim that $\Psi$ is injective. First note that if $E, F \subset M$ are Borel sets and 
\begin{equation}
\label{EF}
m_M(E \cap U_n) = m_M(F \cap U_n), \quad \textrm{for all $n$},
\end{equation}
then $m_M(E \Delta F) = 0$. Indeed, if \eqref{EF} holds, then 
\[
m_M((E \setminus F) \cap U_n) = m_M(E \cap U_n) - m_M(F \cap U_n) = 0
\]
and 
\[
m_M((F \setminus E) \cap U_n) = m_M(F \cap U_n) - m_M(E \cap U_n) = 0,
\]
for all $n$, and thus it suffices to show that if $D \subset M$ is Borel set such
that $m_M(D \cap U_n) = 0$ for all $n$, then $m_M(D)=0$. However, since the 
union of all $U_n$ covers $M$, we must have
\[
m_M(D) = m_M\big(D \cap \big( \bigcup_n U_n \big) \big) \leq \sum_n m_M(D \cap U_n) = 0.
\]
Hence, if $t_1, t_2 \in M$ are such that 
\[
m_M(I t_1 \cap U_n) = m_M(I t_2 \cap U_n), \quad \textrm{for all $n$},
\] 
then $m_M(I \Delta I t_2 t_1^{-1}) = 0$, and thus $t_2 t_1^{-1} \in \Stab_M(I)$, which forces
$t_1 = t_2$, since $\Stab_M(I)$ is assumed to be trivial. This shows that $\Psi$ is injective. \\

Let us now summarize the discussion so far. We have shown that the maps $\Psi$ and $\Phi$ are 
Borel measurable from $M$ and $L$ into $[0,1]^{\bN}$ respectively. Furthermore, $\Psi$ is injective, 
and by \eqref{eqweird}, we have
\[
\Psi(\gamma(s)) = \Phi(s), \quad \forall \, s \in X.
\]
Hence, $\gamma = \Psi^{-1} \circ \Phi$. By Theorem A.4 in \cite{Zi}, $\Psi^{-1}$ 
is Borel measurable, so we conclude that $\gamma$ is Borel measurable as well.
\end{proof}

\subsection{Interlude II: Restrictions of homomorphisms}

The following lemma will be used in \ref{subsubpair} below. 

\begin{lemma}
\label{app}
Suppose that
\begin{itemize}
\item $G$ and $M$ are compact and second countable groups. 
\item $X \subset G$ and $Z \subset X \times X$ are conull.
\item there are Borel measurable maps $\sigma, \tau : X \ra M$ such that
\[
\sigma(x_1) \tau(y_1) = \sigma(x_2) \tau(y_2),
\]
whenever $(x_1,y_1)$ and $(x_2, y_2)$ belong to $Z$ and $x_1 y_1 = x_2 y_2$.
\end{itemize}
Then there exist a continuous homomorphism $\pi : G \ra M$ and $a, b \in M$ such 
that 
\[
a \sigma(x) =\pi(x) \qand \tau(x)b = \pi(x), \quad \textrm{for a.e. $x \in G$}.
\]
\end{lemma}

\begin{proof}
Since $Z \subset X \times X$ is conull, there exist, by Fubini's Theorem, an element $x_o \in X$ and a 
conull subset $X' \subset X$ such that $\{x_o\} \times X' \subset Z$. Since the multiplication map 
$(x,y) \mapsto xy$ pushes $m_G \otimes m_G$ onto $m_G$, we see that the set 
\[
Z' := (x_o,e)^{-1}Z \cap \big\{ (x,y) \in G \times G \, : \, xy \in X' \big\} 
\]
is conull in $G \times G$. We note that for any $(x,y) \in Z'$, we have
\[
(x_o,xy) \in Z \qand (x_ox,y) \in Z,
\]
and thus (since $x_o(xy) = (x_o x)y$),
\begin{equation}
\label{allZ}
\sigma(x_o) \tau(xy) = \sigma(x_o x) \tau(y), \quad \textrm{for all $(x,y) \in Z'$}.
\end{equation}
By Fubini's Theorem, we can find $y_o \in X$ and a conull subset $X'' \subset G$ such that
$X'' \times \{y_o\} \subset Z'$. Arguing as before, we see that the set
\[
Z'' := Z'(e,y_o)^{-1} \cap \big\{ (x,y) \in G \times G \, : \, xy \in X'' \big\}
\]
is conull in $G \times G$. We note that if $(x,y) \in Z''$, then
\[
(x,yy_o) \in Z' \qand (xy,y_o) \in Z',
\]
and thus, by \eqref{allZ}, 
\begin{equation}
\label{wow1}
\sigma(x_o) \tau(xyy_o) = \sigma(x_o x) \tau(yy_o)
\end{equation}
and
\begin{equation}
\label{wow2}
\sigma(x_o) \tau(xyy_o) = \sigma(x_o xy) \tau(y_o),
\end{equation}
for all $(x,y) \in Z''$. \\

Since $X'$ and $X''$ are both conull, so is the set $Y := X'' \cap X' y_o^{-1}$. It is straightforward to check 
that 
$(e,y) \in Z''$ for all $y \in Y$. Hence, by \eqref{wow2}, we have
\[
\sigma(x_o) \tau(yy_o) =\sigma(x_o y) \tau(y_o), \quad \textrm{for all $y \in Y$},
\]
and thus
\begin{equation}
\label{defpio}
\pi_o(y) := \sigma(x_o)^{-1} \sigma(x_o y) = \tau(yy_o) \tau(y_o)^{-1}, \quad \textrm{for all $y \in Y$}.
\end{equation}
Let us now define the conull set
\[
W := \big\{ (x,y) \in Z'' \, : \, xy \in Y \big\} \cap (Y \times Y).
\]
Then, for all $(x,y) \in W$, we have by \eqref{wow1} and \eqref{defpio}, 
\[
\pi_o(xy) = \tau(xyy_o)\tau(y_o)^{-1} = \sigma(x_o)^{-1} \sigma(x_o x) \tau(yy_o) \tau(y_o)^{-1} = \pi_o(x) \pi_o(y).
\] 
In other words, $\pi_o$ satisfies the condition of being a 
homomorphism from $G$ into $M$ almost everywhere. By Theorem B.2 in \cite{Zi}, we can now conclude that there exists a \emph{continuous} homomorphism $\pi : G \ra M$ whose restriction to $W$ coincides with $\pi_o$. In particular, letting $a^{-1} = \sigma(x_o) \pi(x_o)^{-1}$ and $b^{-1} = \pi(y_o)^{-1}\tau(y_o)$, we see from 
\eqref{defpio} that
\[
a \sigma(y) = \pi(y) \qand \tau(y)b = \pi(y), \quad \textrm{for a.e. $y \in G$},
\]
\end{proof}

\subsection{Relations between $\alpha$ and $\beta$}

Let us go back to the setting of Subsection \ref{getrid}, and recall the summary at the end of this subsection. 
In particular, let $S \subset \widehat{N} \setminus \{1\}$ and $\xi \in S$ be as in this subsection. Since
$L$ acts continuously on $N$ by conjugation, it also acts continuously on the dual $\widehat{N}$ via the
adjoint representation,
\[
(s \cdot \eta)(n) = \eta(s^{-1} n s), \quad \textrm{for $\eta \in \widehat{N}$ and $s \in L$}.
\]
In what follows, we will set $\xi_s = s^{-1} \cdot \xi$. This notation should \emph{not} be confused with the notation used in Subsections \ref{subseccomb} and \ref{getrid}. By continuity of the $L$-action on $\widehat{N}$, each
set of the form $\{s \in L \, : \, \xi_s = \eta \big\}$, where $\eta \in \widehat{N}$, is closed in $L$. Hence, since $S$ is 
(at most) countable, the set 
\[
E = \big\{ s \in L \, : \, \xi_s \in S \big\} \subset L
\]
is a countable union of closed subsets of $L$, and thus Borel. In particular, if we define 
the map $\eps : L \ra \{-1,1\}$ by
\[
\eps(s) = 
\left\{
\begin{array}{cc}
+1 & \textrm{if $s \in E$} \\
-1 & \textrm{if $s \notin E$}  
\end{array}
\right.,
\]
then $\eps$ is Borel measurable. Finally, since $\xi \in S$ we have $\eps(e) = 1$. \\

Since we have $\xi^{-1} = \check{\xi}$ for every $\xi \in \widehat{N}$, we note that  $\xi_s^{\eps(s)} \in S$ for all $s \in L$. Furthermore, note that for every Borel set $D \subset \bT$, we have
\begin{equation}
\label{conjxi}
s^{-1}\xi^{-1}(D)s = \xi_s^{-1}(D), \quad \textrm{for all $s \in L$}.
\end{equation}
\subsubsection{Bounding $AB$}

Since $N$ is a normal subgroup of $G$, we have 
\[
sABt^{-1} \cap N = s(ABt^{-1}s \cap s^{-1}Ns)s^{-1} = s(AB(s^{-1}t)^{-1} \cap N)s^{-1}.
\]
for all $s,t \in L$. Hence, by \eqref{summaryAB} and \eqref{conjxi},
\[
AB(s^{-1}t)^{-1} \cap N \sim s^{-1}\xi^{-1}(IJ \alpha(s) \beta(t))s = \xi_s^{-1}(IJ \alpha(s) \beta(t)),
\]
for all $(s,t) \in Y_o$. We see that the left hand side only depends on $s^{-1}t$. In particular, for all pairs $(s,t)$ and $(u,v)$ in $Y_o$ such that $s^{-1}t = u^{-1}v$, 
we must have
\[
\xi_s^{-1}(IJ  \alpha(s) \beta(t))) \sim \xi_u^{-1}(IJ\alpha(u) \beta(v)).
\]
Since $N$ is normal, we conclude from this that 
\[
\ker \xi_s = \ker \xi_u = \ker \xi,
\]
and thus the composition $\xi_s \circ \xi_u^{-1}$ is a well-defined automorphism of $\bT$. Since $\Aut(\bT) = \{-1,1\}$,
we see that \emph{either}
\[
\xi_s = \xi_u \qand \alpha(s) \beta(t) = \alpha(u) \beta(v)
\]
or
\[
\xi_s = \check{\xi}_u \qand \alpha(s) \beta(t) = (\alpha(u) \beta(v))^{-1}.
\] 
Since $\xi_s^{\eps(s)} \in S$ for all $s \in L$, we conclude that 
\begin{equation}
\label{sueq}
\xi_s^{\eps(s)} = \xi_u^{\eps(u)} 
\qand
(\alpha(s) \beta(t))^{\eps(s)} = (\alpha(u) \beta(v))^{\eps(u)},
\end{equation}
whenever $(s,t), (u,v) \in Y_o$ with $s^{-1}t = u^{-1}v$. Since the first identity is independent of $t$ and $v$, 
we see that the map $s \mapsto \xi_s^{\eps(s)}$ is almost everywhere constant, say equal to $\eta \in \widehat{N} \setminus \{1\}$ on a conull subset of $X_o$ (which we henceforth identify with $X_o$). We note that the sets
\[
R_{+} = \big\{ s \in L \, : \, \xi_s = \eta \big\} \qand R_{-} = \big\{ s \in L \, : \, \xi_s^{-1} = \eta \big\}
\]
are \emph{closed}, and by assumption $X_o \subset R_{+} \cup R_{-}$. Since $X_o$ is conull, and the union $R_{+} \cup R_{-}$ is closed, we conclude that $R_{+} \cup R_{-} = L$. Hence, $\xi_s^{\eps(s)} = \eta$ for \emph{all} $s \in L$.
Since $\eps(e) = 1$, we see that $\xi = \eta$, and thus
\begin{equation}
\label{whatepsreallydoes}
\xi(sns^{-1})^{\eps(s)} = \xi(n), \quad \textrm{for all $s \in L$ and $n \in N$}.
\end{equation}
We conclude that the $L$-action on $\widehat{N}$ preserves the two-element set $\{ \xi, \check{\xi} \}$, and the
corresponding homomorphism $L \ra \Aut(\bT) \cong \{-1,1\}$, coincides with $\eps$. Since $\eps$ is a Borel measurable
homomorphism between second countable groups, it must be continuous by Theorem B.3 \cite{Zi}. \\
%
%
%
%

\emph{To summarize:} 
There exists a continuous homomorphism $\eps : L \ra \{-1,1\}$ such that the relation \eqref{whatepsreallydoes} holds and
\begin{equation}
\label{abeps}
(\alpha(s) \beta(t))^{\eps(s)} = (\alpha(u) \beta(v))^{\eps(u)},
\end{equation}
whenever $s^{-1}t = u^{-1}v$. In particular, for any Borel set $D \subset \bT$, we have
\begin{equation}
\label{whatepsdoes}
s^{-1} \xi^{-1}(D) s = \xi_s^{-1}(D) 
= 
\xi^{-1}(D^{\eps(s)}),
\end{equation}
for all $s \in L$, where we adopt the convention that $D^1 = D$. 

\subsubsection{Bounding $A$}

By \eqref{summary}, we have
\[
sA \cap N = s(A \cap s^{-1}N) \subset \xi^{-1}(I \alpha(s)), \quad \textrm{for all $s \in X_o$},
\]
and thus, by \eqref{whatepsdoes} and our assumption that $I$ is symmetric, 
\[
A \cap Ns^{-1} \subset \big(s^{-1}\xi^{-1}(I \alpha(s))s \big) s^{-1} = \xi^{-1}(I \alpha(s)^{\eps(s)})s^{-1},
\]
for all $s \in X_o$. Since $X_o^{-1} = X_o$ and $\eps(s^{-1}) = \eps(s)$ (indeed, $1 = \eps(ss^{-1}) = \eps(s) \eps(s^{-1})$ for all $s$), we also have
\begin{equation}
\label{ANSIes}
A \cap Ns \subset \xi^{-1}(I \alpha(s^{-1})^{\eps(s)})s, \quad \textrm{for all $s \in X_o$}.
\end{equation}
Let us define the map $\sigma : NX_o \ra \bT \rtimes \{-1,1\}$ by
\begin{equation}
\label{defsigma}
\sigma(ms) = \big(\xi(m) \alpha(s^{-1})^{-\eps(s)},\eps(s)).
\end{equation}
Since both $\alpha$ and $\eps$ are Borel measurable, so is $\sigma$. We note that
\[
\sigma^{-1}(I \rtimes \{-1,1\}) \cap Ns= \big\{ ms \, : \, \xi(m) \alpha(s^{-1})^{-\eps(s)} \in I \big\} = \xi^{-1}(I \alpha(s^{-1})^{\eps(s)})s,
\]
for all $s \in X_o$, and thus, by using this and \eqref{ANSIes}, we get 
\[
A \cap Ns \subset \sigma^{-1}(I \rtimes \{-1,1\}) \cap Ns, \quad \textrm{for all $s \in X_o$}.
\]
We conclude that
\begin{equation}
\label{defA'}
A' := A \cap NX_o \subset \sigma^{-1}(I \rtimes \{-1,1\}).
\end{equation}
Note that $A'$ is a conull subset of $A$, since $NX_o$ is a conull subset of $G$.

\subsubsection{Bounding $B$}

By \eqref{summary}, we have
\[
Bt^{-1} \cap N = (B \cap Nt)t^{-1} \subset \xi^{-1}(J\beta(t)),
\]
for all $t \in X_o$, and thus
\[
B \cap Nt \subset \xi^{-1}(J\beta(t))t, \quad \textrm{for all $t \in X_o$}.
\]
Define the map $\tau : NX_o \ra \bT \rtimes \{-1,1\}$ by
\begin{equation}
\label{deftau}
\tau(nt) = (\xi(n)\beta(t)^{-1},\eps(t)).
\end{equation}
Since both $\beta$ and $\eps$ are Borel measurable, so is $\tau$. Note that
\[
\tau^{-1}(J \rtimes \{-1,1\}) \cap Nt = \big\{ nt \, : \, \xi(n) \beta(t)^{-1} \in J \big\} = \xi^{-1}(J\beta(t))t,
\]
for all $t \in X_o$, and thus
\[
B \cap Nt \subset \tau^{-1}(J \rtimes \{-1,1\}) \cap Nt, \quad \textrm{for all $t \in X_o$}.
\]
We conclude that
\begin{equation}
\label{defB'}
B' := B \cap NX_o \subset \tau^{-1}(J \rtimes \{-1,1\}),
\end{equation}
and $B' \subset B$ is conull, since $NX_o$ is a conull subset of $G$.

\subsubsection{The pair $(\tau,\sigma)$}
\label{subsubpair}
We wish to verify the pair $(\sigma, \tau)$ of Borel maps above satisfies the conditions in Lemma \ref{app} with $M = \bT \rtimes \{-1,1\}$ and the conull subsets
\[
X = NX_o \subset G \qand Z = \big\{ (ms,nt) \, : \, m,n \in N, \enskip (s,t) \in Y_o \big\} \subset X \times X.
\] 
The multiplication in $\bT \rtimes \{-1,1\}$ of two elements $(r_1,\delta_1)$ and $(r_2,\delta_2)$ will be written
\begin{equation}
\label{multdef}
(r_1,\delta_1)(r_2,\delta_2) = (r_1 r_2^{\delta_1},\delta_1 \delta_2).
\end{equation}
Suppose that $(ms,nt)$ and $(pu,qv)$ belong to $Z$ and
\[
(ms)(nt) = m(sns^{-1})st = (pu)(qv) = p(uqu^{-1})uv.
\]
Since $N \cap L = \{e\}$, this forces
\[
m(sns^{-1}) = p(uqu^{-1}) \qand st = uv.
\]
By \eqref{abeps} and our assumption that $X_o^{-1} = X_o$, we have
\[
(\alpha(s^{-1}) \beta(t))^{\eps(s)} = (\alpha(u^{-1}) \beta(v))^{\eps(u)},
\]
whenever $(s,t)$ and $(u,v)$ belong to $Y_o$ and $st = uv$. \\

Recall \eqref{whatepsreallydoes} and the definitions of $\sigma$ and $\tau$ from \eqref{defsigma} and \eqref{deftau} respectively. Upon combining the relations above, and using the multiplication convention in $\bT \rtimes \{-1,1\}$ 
explained in \eqref{multdef}, we see that
\begin{eqnarray*}
\sigma(ms) \tau(nt) 
&=&
\big(\xi(m) \alpha(s^{-1})^{-\eps(s)},\eps(s)) \, (\xi(n)\beta(t)^{-1},\eps(t)) \\
&=&
\big(\xi(m) \xi(n)^{\eps(s)}\alpha(s^{-1})^{-\eps(s)}\beta(t)^{-\eps(s)}, \eps(s) \eps(t)\big) \\
&=&
\big(\xi(m) \xi(sns^{-1}) (\alpha(s^{-1})\beta(t))^{-\eps(s)}, \eps(st) \big) \\
&=&
\big(\xi(msns^{-1}) (\alpha(s^{-1})\beta(t))^{-\eps(s)}, \eps(st) \big) \\
&=&
\big(\xi(puqu^{-1}) (\alpha(u^{-1})\beta(v))^{-\eps(u)}, \eps(uv) \big) \\
&=&
\sigma(pu) \tau(qv). 
\end{eqnarray*}
By Lemma \ref{app} we conclude that there exist a continuous homomorphism $\pi : G \ra \bT \rtimes \{-1,1\}$
and $a, b \in M$ such that $\sigma(g) = a^{-1}\pi(g)$ and $\tau(g) = \pi(g) b^{-1}$ almost everywhere. Upon 
possibly passing to further conull subsets in \eqref{defA'} and \eqref{defB'}, we conclude that
\begin{equation}
\label{defAB'}
A' \subset \pi^{-1}(a(I \rtimes \{-1,1\})) \qand B' \subset \pi^{-1}((I \rtimes \{-1,1\})b),
\end{equation}
where $A' \subset A$ and $B' \subset B$ are conull subsets. 

\subsubsection{Determining possible images of $\pi$}
\label{deter}
We recall that $M = \bT \rtimes \{-1,1\}$ and $\pi$ is a continuous homomorphism of $G = N \rtimes L$ into $M$. 
Since $N$ is connected and abelian, $\pi(N)$ is a compact and connected abelian subgroup of $M$, and thus 
either trivial or equal to $\bT$. We claim that the first case cannot occur. Indeed, recall that our standing 
assumption in Proposition \ref{relativemain} is that $(A,B)$ is sub-critical with respect to $N$, which by 
\eqref{subAandBL} in particular implies that 
\[
m_N(sA \cap N) = m_G(A), \quad \textrm{for $m_L$-a.e. $s \in L$}.
\]
Since $A' \subset A$ is conull, and
\[
\int_L m_N(sA' \cap N) \, dm_L(s) = m_G(A') = m_G(A) = \int_L m_N(sA' \cap N) \, dm_L(s),
\]
we see that $ m_N(sA' \cap N) = m_G(A)$ for $m_L$-a.e. $s \in L$ as well. \\

By \eqref{defAB'}, the set $A'$
is also a conull subset of $\pi^{-1}(a(I \rtimes \{-1,1\})$, so the same type of argument as above shows that 
\[
m_N(s \pi^{-1}(a(I \rtimes \{-1,1\}) \cap N)) = m_G(A), \quad \textrm{for $m_L$-a.e. $s \in L$}.
\]
We now note that if $\pi(N) = \{e_M\}$, so that $N < \ker \pi$, then the left-hand side is either $0$ or
$1$, which contradicts our assumption that $0 < m_G(A) < 1$. 

We conclude that $\pi(N) = \bT$, so we have two possibilities: Either $\pi(L) \supset \{0\} \rtimes \{-1,1\}$, in which
case we must have $\pi(G) = \bT \rtimes \{-1,1\}$, or $\pi(G) = \bT \rtimes \{1\} \cong \bT$. 
 
\subsubsection{Finishing the proof of Proposition \ref{relativemain}}

Let us briefly summarize the argument so far: We have produced conull subsets $A' \subset A$ and $B' \subset B$
such that 
\[
A' \subset \pi^{-1}(a(I \rtimes \{-1,1\})) \qand B' \subset \pi^{-1}((I \rtimes \{-1,1\})b),
\]
where $I, J \subset \bT$ are closed intervals with
\[
m_G(A) = m_{\bT}(I) \qand m_G(B) = m_{\bT}(J).
\]
From the previous subsection, we know that $\pi(G)$ can be either $\bT$ or $\bT \rtimes \{-1,1\}$, and thus in 
either case, 
\[
m_G(\pi^{-1}(a(I \rtimes \{-1,1\})) = m_\bT(I) \qand m_G(\pi^{-1}((J \rtimes \{-1,1\})b) = m_\bT(J).
\]
We want to prove that $(A',B')$ reduces to a Sturmian pair in either $\bT$ or $\bT \rtimes \{-1,1\}$. From
the arguments above, 
it is clear that it remains to show that $m_G(A'B') = m_G(AB)$. We shall argue by contradiction: If $m_G(A'B') < m_G(AB)$, then, since $A' \subset A$ and $B' \subset B$ are conull subsets, we have
\[
m_G(A' B') < m_G(AB) = m_G(A) + m_G(B) = m_G(A') + m_G(B') < 1,
\]
and thus $(A',B')$ is critical in $G$. By the first assertion 
in Theorem \ref{knes}, this would imply that there exists a finite factor group of $G$ such that neither $A'$ nor $B'$ project onto it. This contradicts our assumption that $(A,B)$ is spread-out in $G$, which finishes the proof of Proposition \ref{relativemain}.

\section{Acknowledgements}

The author would like to thank the referee for a very thorough read of the present paper, as well as previous,
less accurate, versions thereof.  

\appendix

\section{Compact groups with a dense amenable subgroup}

A countable group $\Gamma$ is \textbf{amenable} if there is sequence $(F_n)$ of finite subsets of $\Gamma$ 
such that
\[
\varlimsup_n \frac{|F_n \Delta \gamma F_n|}{|F_n|} = 0, \quad \textrm{for all $\gamma \in \Gamma$}.
\]
It is well-known, see e.g. the book \cite{Pat}, that every countable solvable (so in particular every abelian or nilpotent) 
group is amenable. Furthermore, subgroups and quotients of amenable groups are again amenable. On the other hand, countable groups with a free subgroup of rank at least two cannot be 
amenable. It was shown by Tits \cite{Tits}, that the absence of a free subgroup of rank two exactly characterizes
amenability among linear groups. \\

In the paper \cite{BF6} by the author and A. Fish, "density analogues" of the results by Kemperman and Kneser
mentioned in the introduction are established for general countable amenable groups. This is done by a technical reduction of product sets in a given countable amenable group $\Gamma$ to product sets in an associated 
(metrizable) compactification of $\Gamma$, i.e. a compact and metrizable group $K$ which contain $\Gamma$
as a dense subgroup. As it turns out, the condition to contain a dense countable amenable group, puts serious
constraints on the identity component of $K$ - It has to be abelian. This observation motivated the study pursued 
in this paper.  For completeness, we sketch a proof here. 

\begin{proposition}
\label{denseamen}
If $K$ is a compact and second countable Hausdorff group with a dense countable amenable subgroup $\Gamma$, then the identity component $K^o$ of $K$ is abelian.
\end{proposition}

\begin{proof}
We shall argue by contradiction: Suppose that there exist two elements $x,y \in K^o$ such that $xyx^{-1}y^{-1} \neq e_K$. By Peter-Weyl's Theorem, we can find a positive integer $n$ and a representation $\pi$ of $K$ into $U(n)$ such that $\pi(x)\pi(y)\pi(x)^{-1}\pi(y)^{-1} \neq e_{\pi(K)}$. We note that $\pi(x)$ and $\pi(y)$ belong to
the identity component of the (possibly not connected) compact Lie group $\pi(K)$. Let $\Gamma_\pi$ denote the image of 
$\Gamma$ under $\pi$; by assumption $\Gamma_\pi$ is a dense countable amenable subgroup of $\pi(K)$.
By Theorem 6.5 (iii) in \cite{HM}, $\pi(K)^o$ has finite index in $\pi(K)$, and a straightforward argument shows that
$\Lambda := \Gamma_\pi \cap \pi(K)^o$ is a dense countable amenable subgroup of $\pi(K)^o$. In particular, 
the commutator subgroup $[\Lambda,\Lambda]$ is a dense amenable subgroup of $[\pi(K)^o,\pi(K)^o]$. By
Theorem 6.18 in \cite{HM}, the latter group is a semisimple and connected compact Lie group. At this point, 
Tits' Alternative \cite{Tits} can be applied: a non-trivial semisimple and connected compact Lie group cannot contain a dense amenable subgroup. 
We conclude that $\pi(K)^o$ is abelian, which contradicts $\pi(x)\pi(y)\pi(x)^{-1}\pi(y)^{-1} \neq e_{\pi(K)}$.
\end{proof}

\end{document}